\documentclass{bcp}

\def\LaTeX{L\kern -.36em\raise .3ex\hbox{\sc a}\kern -.15em T\kern -.1667em%
\lower .7ex\hbox{E}\kern -.125em X}

\usepackage{amssymb}
\usepackage{amsmath}
\usepackage{latexsym}
\usepackage[dvips]{graphics}
\usepackage{amscd}
\usepackage{amsfonts}
\usepackage{theorem}
\usepackage[all]{xy}
\usepackage{mdwtab}
\usepackage{float}

{\theorembodyfont{\rmfamily} 

}

\newcommand{\RR}{\mathbb R}
\newcommand{\RP}{\mathbb{RP}}
\def\C{\mathrm {C}}
\def\N{\mathrm {N}}
\def\K{\mathrm {K}}
\def\H{\mathrm {H}}
\def\MC{\mathrm {MC}}
\def\MCC{\mathrm {MCC}}
\newcommand{\eps}{\varepsilon}
\renewcommand{\d}{\partial}
\renewcommand{\:}{\colon\,}
\newcommand{\id}{\mathrm {id}}
\newcommand{\pinch}{\mathrm {pinch}}
\newcommand{\Int}[1]{{\stackrel{{\scriptscriptstyle\circ}}{#1}}}
\newcommand{\qed}{\mbox{}\hfill$\square$}

\newcommand{\Bigdownarrow}{\Downarrow}

\begin{document}

 \keywords{Coincidence invariants; Nielsen numbers; Kervaire invariants}
 \mathclass{Primary 54H25, 55M20; Secondary 55P35, 55Q15, 55Q22, 55Q40.}
 \abbrevauthors{U.\ Koschorke}
 \abbrevtitle{Some homotopy theoretical questions}

\title{Some homotopy theoretical questions \\ arising in Nielsen coincidence theory}

\author{Ulrich Koschorke}
\address{Fachbereich Mathematik, Universit\"at Siegen\\
57068 Siegen, Germany\\
E-mail: koschorke@mathematik.uni-siegen.de}

\maketitlebcp

\abstract{Basic examples show that coincidence theory is intimately
related to central subjects of differential topology and homotopy
theory such as Kervaire invariants and divisibility properties of
Whitehead products and of Hopf invariants. We recall some recent
results and ask a few questions which seems to be important for a
more comprehensive understanding.}

\bigskip
\bigskip

Throughout this paper $M$ and $N$ will denote closed connected smooth
manifolds of dimensions $m$ and $n\ge2$, resp.

\defin{Definition}{1 (cf.~\cite {K2}, (2), (3) and~1.1).
Given (continuous) maps $f_1,f_2\:M\to N$, let $\MC(f_1,f_2)$ and
$\MCC(f_1,f_2)$, resp., denote the {\it minimum number}\/ of points
and of path components, resp., among all the coincidence subspaces
$$
\C(f_1',f_2') = \{x\in M\mid f_1'(x)=f_2'(x)\} \ \subset \ M
$$
of maps $f_i'$ which are homotopic to $f_i$, $i=1,2$.

The pair $(f_1,f_2)$ is called {\it loose}\/ if $\MC(f_1,f_2)=0$
(or, equivalently, $\MCC(f_1,f_2)=0$), i.e.\ if the maps $f_1$ and
$f_2$ can be deformed away from one another.}

The principal aim in topological coincidence theory is to get a good
understanding of the minimum numbers $\MC(f_1,f_2)$ and
$\MCC(f_1,f_2)$ (compare~\cite {B}, p.~9) and, in particular, to
obtain precise looseness criteria.

In~\cite {K2} a geometric invariant $\omega^\#(f_1,f_2)$ was
introduced: it is the bordism class of the (generic) coincidence
submanifold
$$ \refstepcounter {Thm} \label {eq:2}
C=\C(f_1,f_2) = \{x\in M\mid f_1(x)=f_2(x)\}
 \leqno (\theThm)
$$
of $M$, together with a description of its normal bundle and a map
into a certain path space $E(f_1,f_2)$. Moreover, a simple numerical
invariant was extracted: the {\it Nielsen number} $\N^\#(f_1,f_2)$,
which counts the ``essential'' pathcomponents of the target space
$E(f_1,f_2)$. We have
$$ \refstepcounter {Thm} \label {eq:3}
\N^\#(f_1,f_2)\le \MCC(f_1,f_2)\le \MC(f_1,f_2);
 \leqno (\theThm)
$$
in many cases, e.g.\ in the {\it stable range} $m\le2n-3$, actually
$\N^\#(f_1,f_2)=\MCC(f_1,f_2)$ (compare~\cite{K1}, 1.10).

\defin{Example}{2 (Fixed point theory).
Here $M=N$, $f_2=$ identity map, and $\N^\#(f_1,\id)$ is the
classical Nielsen number of $f_1$; it is known to coincide with the
minimum number of fixed points occuring in the homotopy class of
$f_1$, provided $n\ne2$ or $\chi(N)\ge0$ (``Wecken theorems'',
cf.~\cite {B}). \qed}

In this paper we will discuss the case $M=S^m$ which is particulary
accessible to the methods of homotopy theory. Here we can identify
$\omega^\#(f_1,f_2)$ with the bordism class of the (generic)
coincidence submanifold $\C(f_1,f_2)$ of $S^m$ together with a
(nonstabilized) framing and a map into the loop space $\Omega N$ of
$N$. Via the Pontryagin--Thom procedure, this can be interpreted by
maps from $S^m$ into the Thom space of the trivial $n$-plane bundle
over $\Omega N$. We obtain the homomorphism
$$ \refstepcounter {Thm} \label {eq:5}
\omega^\#\:\pi_m(N) \oplus \pi_m(N) \longrightarrow
\pi_m(S^n\wedge(\Omega N)^+)
 \leqno (\theThm)
$$
where $(\Omega N)^+$ denotes the loop space with one point added
disjointly (cf.\ section~6 of~\cite {K2}; for a stabilized version
and its relations e.g.\ to Hopf-James invariants see~\cite{K1},
1.14).

From now on we consider only maps $f_1,f_2,f\ldots\:S^m\to N$ where
we assume that the order
 $$
 k:=\#\pi_1(N)
 $$
of the fundamental group of $N$ is finite (otherwise $(f_1,f_2)$ is
loose for all $f_1,f_1\:S^m\to N$, cf.~\cite{K4}, 1.3). Then the
exact homotopy sequence of the obvious projection
$\bigvee^kS^n\vee\tilde N\to\tilde N$, turned into a fibre map,
yields the isomorphism
$$ \refstepcounter {Thm} \label {eq:6}
\kappa \: \pi_m(S^n\wedge(\Omega N)^+) \longrightarrow \pi_m
(\bigvee^kS^n\vee\tilde N,\tilde N)
 \leqno (\theThm)
$$
(cf.~\cite {K2}, (61)) which may give useful information on the
target group of $\omega^\#$; here $\tilde N$ denotes the universal
covering space of $N$.

If $N$ happens to be a spherical space form (so that $\tilde N\cong
S^n$), we can describe the homotopy groups of the wedge
$\bigvee^kS^n\vee\tilde N$ as a direct sum of the homotopy groups of
spheres $S^{j(n-1)+1}$ (cf.~\cite {H}); often this allows explicit
calculations.

\th{Question}{A.}{What can be said about the homotopy groups of the
wedge $\bigvee^kS^n\vee\tilde N$ for general manifolds?}

We have the logical implications
$$ \refstepcounter {Thm} \label {eq:7}
(f_1,f_2) \mbox{ is loose} \ \implies \ \omega^\#(f_1,f_2)=0 \ \iff \
\N^\#(f_1,f_2)=0
 \leqno (\theThm)
$$
(cf.~\cite {K4}, theorem 1.30). Of course it is of central importance
to know when $\omega^\#$ and $\N^\#$ are {\it complete} looseness
obstructions, i.e.\ when the vanishing of \ $\omega^\#(f_1,f_2)$ or,
equivalently, of $\N^\#(f_1,f_2)$ is also sufficient for the pair
$(f_1,f_2)$ to be loose. Since
$$ \refstepcounter {Thm} \label {eq:8}
\omega^\#([f_1],[f_2])\ \ =\ \ \omega^\#([f_1]-[f_2],0)\ +\
\omega^\#([f_2],[f_2])
 \leqno (\theThm)
$$
the following two settings are of particular interest.

\medskip
{\bf I. The root case: $f_2\equiv y_0$ } (where the fixed value
$y_0\in N$ is given).

Here our invariant yields the degree homomorphism
$$ \refstepcounter {Thm} \label {eq:9}
\deg^\#:=\omega^\#(-,y_0)\ \:\ \pi_m(N)\to\pi_m(S^n\wedge(\Omega
N)^+)
 \leqno (\theThm)
$$
which turns out to be essentially an enriched Hopf--Ganea--invariant
homomorphism (see~\cite {K2}, theorem~7.2). It is also related to the
homomorphism
$$ \refstepcounter {Thm} \label {eq:10}
(\pinch_*,\d)\: \pi_m(\tilde N,\tilde N\setminus\bigcup^k \Int{B^n})
\longrightarrow \pi_m(\bigvee^kS^n\vee\tilde N, \tilde N) \oplus
\pi_{m-1}(\tilde N\setminus \bigcup^k \Int{B^n})
 \leqno (\theThm)
$$
(cf.~\cite {K2}, 7.3) where $\bigcup^k \Int{B^n}$ is the disjoint
union of open $n$-balls whose (embedded) boundary spheres intersect
only in $\tilde y_0\in\tilde N$; the ``pinching map'' pinch collapses
these boundary spheres into a single point and deforms $\tilde
N\setminus \bigcup B^n$ to $\tilde N\setminus\{\tilde y_0\}$; $\d$
denotes the obvious connecting homomorphism.

Clearly $\deg^\#$ vanishes on $i_*(\pi_m(N\setminus\{y_0\}))$ where
$i$ denotes the inclusion $N\setminus\{y_0\}\subset N$.

The following conditions are equivalent (cf.~\cite {K2}, 6.4
and~7.3):

(i) the sequence
$$
\pi_m(N\setminus\{y_0\})
 \stackrel {i_*}{\longrightarrow} \pi_m(N)
 \stackrel {\deg^\#}{\longrightarrow} \pi_m(S^n\wedge(\Omega N)^+)
$$
is exact;

(ii) the homomorphism $(\pinch,\d)$ (cf.~(\ref {eq:10})) is
injective;

(iii) $\deg^\#(f)$ is the {\it complete} looseness obstruction for
all pairs of the form $(f,y_0)$, where $f\:S^m\to N$.

All these conditions are often satisfied, e.g.\ when $m\le 2n-3$ or
$n=2$ or $N$ is a sphere or a (real, complex or quaternionic)
projective space of arbitrary dimension (cf.~\cite {K2},
theorem~6.5). Do they always hold?

\th{Question}{B.}{What can be said about the kernel of $(\pinch,\d)$
(cf.~(\ref {eq:10})) for general $N$ and arbitrary dimensions?}

It may also be interesting to note that the Nielsen number
$\N^\#(f,y_0)$ can assume only the values $0$ and $k$ (cf.~\cite
{K2}, 4.3).

\medskip
{\bf II. The selfcoincidence case: $f_1=f_2=:f$. }

Here the loopspace aspect of our invariants carries no extra
information. Therefore $\omega^\#(f,f)$ is precisely as strong as
its image in $\pi_m(S^n)$ under the obvious forgetful map, and the
Nielsen Number $\N^\#(f,f)$ can take only the values $0$ and $1$.
However, we can distinguish between two kinds of deformations: small
deformations (which move $f$ only $\eps$-far away for a small
$\eps>0$) on one hand, and arbitrary deformations (which may use all
the space available in $N$) on the other hand (cf.~\cite {DG}).

Consider the bundle $ST(N)$ of unit tangent vectors over $N$ (with
respect to some metric) and the resulting exact (horizontal) homotopy
sequence as well as the Freudenthal suspension homomorphism $E$:
$$ \refstepcounter {Thm} \label {eq:11}
\xymatrix{
 \ldots \ar[r] & \pi_m(ST(N)) \ar[r] & \pi_m(N) \ar[r]^{\!\!\!\!\!\!\d} & \pi_{m-1}(S^{n-1}) \ar[r] \ar[d]^E & \ldots
 \\
               &                     &                    & \pi_m(S^n)                         & .
}
 \leqno (\theThm)
$$

\th{Theorem}{3.}{Given $[f]\in\pi_m(N)$, we have the following
logical implications:

\medskip
{\rm (i)} $\d([f])\in\pi_{m-1}(S^{n-1})$ vanishes;

$\phantom{j}\!$$\Updownarrow$

{\rm (ii)} $(f,f)$ is loose by small deformation;

\medskip
$\phantom{i}$$\Bigdownarrow$ \quad $\left(\Updownarrow \mbox { if }
N=\RP(n)\right)$

\medskip
{\rm (iii)} $(f,f)$ is loose (by any deformation);

\medskip
$\phantom{i}$$\Bigdownarrow$ \quad $\left(\Updownarrow \mbox { if }
N=S^n\right)$

\medskip
{\rm (iv)} $\omega^\#(f,f)=0$;

$\phantom{i}$$\Updownarrow$

{\rm (v)} $E(\d([f])=0$.}

Thus $\omega^\#(f,f)$ is just ``one desuspension short'' of being the
{\it complete} looseness obstruction.

The equivalence of~(i) and~(ii) was already noted by A.~Dold and
D.L.~Gon\c calves in~\cite {DG}.

Observe also that all the conditions (i)$,\ldots,$(v) above {\it
except} (iii) are compatible with covering projections $p\:\tilde
N\to N$ (compare~\cite {K3},1.22).

\th{Corollary}{4.}{The conditions {\rm (i)}$,\ldots,${\rm (v)} in
theorem~$3$ are all equivalent if the suspension homomorphism $E$,
when restricted to $\d(\pi_m(N))$ (cf.~$(\ref {eq:11})$), is
injective and, in particular, if $m\le n+3$ or if $m=n+4\ne10$ or in
the stable dimensional range $m\le 2n-3$.}

Indeed, in these three dimension settings $E$ is injective whenever
$n\equiv 0(2)$ (compare e.g.~\cite {T} or~\cite {K4}, 4.5). If the
Euler characteristic of $N$ vanishes (e.g.\ when $n\not\equiv0(2)$),
then the conditions (i)$,\ldots,$(v) in theorem~3 are automatically
satisfied due to the existence of a nowhere vanishing vector field
along which the map $f$ can be pushed slightly away from itself.

\th{Corollary}{5.}{Assume that

{\rm (i)} $\#\pi_1(N)>2$ and $N$ is orientable or not; or

{\rm (ii)} $\#\pi_1(N)\ge2$ and $N$ is orientable.

If in addition $E|_{\d(\pi_m(N))}$ is injective, then for all maps
$f\:S^m\to N$ the pair $(f,f)$ is loose by small (and hence by any)
deformation.}

Indeed, according to theorem~1.21 of~\cite {K4} \ $\omega^\#(f,f)$
vanishes here.

This suggests that our invariants are particularly interesting when
$N$ has a small fundamental group, e.g.\ in the case of spheres and
projective spaces.

\th{Theorem}{6.}{Given a map $f\:S^m\to\RP(n)$, $m,n\ge 2$, let
$\tilde f\:S^m\to S^n$ be a lifting. Then the following conditions
are equivalent:

{\rm (i)} $\omega^\#(f,f)=0$, but $(f,f)$ is not loose;

{\rm (ii)} $(\tilde f,\tilde f)$ is loose but not by small
deformation;

{\rm (iii)} $\d_\RR([\tilde f])\ne0$, but $E\circ\d_\RR([\tilde
f])=0$;

{\rm (iv)} $(\tilde f,\tilde f)$ is loose, but $(f,f)$ is not loose;

{\rm (v)} $\MC(\tilde f,\tilde f)<\MC(f,f)$;

{\rm (vi)} $\MCC(\tilde f,\tilde f)<\MCC(f,f)$.

In particular, maps $f,\tilde f$ which satisfy these conditions exist
precisely in those dimension combinations where $E$ is not injective
on $\d(\pi_m(S^n))$.}

Already in the first nonstable dimension setting we encounter
fascinating interrelations with other, seemingly distant, branches
of topology.

\th{Theorem}{7.}{Let $\tilde f\:S^{2n-2}\to S^n$ be a lifting of a
map $f\:S^{2n-2}\to \RP(n)$. Assume that $n$ is even, $n\ne 2,4,8$.

Then the pair $(f,f)$ is loose if and only if both \
$\omega^\#(f,f)$ and the Kervaire invariant $\K([\tilde f])$ of
$\tilde f$ vanish.}

Such a connection with the Kervaire invariant was originally
discovered by D.L.~Gon\-\c cal\-ves and D.~Randall~\cite {GR2} who
also studied the second nonstable dimension setting (in~\cite {GR1})
and found the first examples illustrating a version of the following
result (cf.~\cite {KR}).

\th{Theorem}{8.}{Let $\tilde f\:S^{2n-1}\to S^n$ be a lifting of a
map $f\:S^{2n-1}\to\RP(n)$. Assume that $n\equiv2(4)$, $n\ge6$.

Then $(f,f)$ is loose if and only if $\omega^\#(f,f)=0$ and, in
addition, the Hopf invariant $\H(\tilde f)$ is divisible by $4$.}

In the next nonstable dimension settings the noninjectivity of
$E|_{\d(\pi_m(S^n))}$ turns out to be closely related to the
question whether the Whitehead products of $\iota_{n-1}$ with
$\eta^2_{n-1},\nu_{n-1},\nu^2_{n-1},\sigma_{n-1},\dots$ can be
halved, i.e.\ lie in $2\pi_*(S^{n-1})$. Many relevant results have
been listed by M.~Golasi\'nski and J.~Mukai in their very helpful
paper~\cite {GM}, but what else is known?

\th{Question}{C.}{Which Whitehead products in $\pi_{m-1}(S^{n-1})$,
$n$ even, are divisible by $2$?}

This seems to be a subtle problem; in coincidence theory it arises
already in the very special case where $N$ equals $S^n$ or $\RP(n)$.

\th{Question}{D.}{What can be said about the subgroup $\d(\pi_m(N))$
of $\pi_{m-1}(S^{n-1})$ (cf.~(\ref {eq:11})) for arbitrary manifolds
$N$ with nontrivial Euler characteristic?

When is $E|_{\d(\pi_m(N))}$ injective?}

E.g.\ what about the case where $N$ is a general Grassmanian
manifold?

\end{document}